\documentclass{amsart}
\usepackage{here,graphicx,color,amsthm,amsmath,amssymb}
\usepackage[frame,cmtip,curve,arrow,matrix,line,graph]{xy}
\newcommand{\nc}{\newcommand}

\newtheorem{thm}{Theorem}[section]
\newtheorem{rmk}[thm]{Remark}
\newtheorem{prop}[thm]{Proposition}
\newtheorem{lemma}[thm]{Lemma}

\newtheorem{corollary}[thm]{Corollary}

\newtheorem{definition}{Definition}[section]

\newenvironment{defin}{\begin{definition} \rm}{\end{definition}}
\newenvironment{cor}{\begin{corollary} \rm}{\end{corollary}}
\newenvironment{lem}{\begin{lemma}\rm }{\end{lemma}}
\nc{\Ext}{\operatorname{Ext}}
\nc{\Rep}{\operatorname{Rep}}
\nc{\Hom}{\operatorname{Hom}}
\nc{\Coh}{\operatorname{Coh}}
\nc{\Aut}{\operatorname{Aut}}
\nc{\D}{\operatorname{D}}
\nc{\Stab}{\operatorname{Stab}}
\nc{\GL}{\operatorname{GL}}
\nc{\Log}{\mathop{\mathrm{Log}}}
\nc{\abs}[1]{\lvert#1\rvert}

\nc{\cA}{{\mathcal A}}
\nc{\cB}{{\mathcal B}}
\nc{\cC}{{\mathcal C}}
\nc{\cD}{{\mathcal D}}
\nc{\cE}{{\mathcal E}}
\nc{\cF}{{\mathcal F}}
\nc{\cG}{{\mathcal G}}
\nc{\cH}{{\mathcal H}}
\nc{\cI}{{\mathcal I}}
\nc{\cJ}{{\mathcal J}}
\nc{\cK}{{\mathcal K}}
\nc{\cL}{{\mathcal L}}
\nc{\cM}{{\mathcal M}}
\nc{\cN}{{\mathcal N}}
\nc{\cO}{{\mathcal O}}
\nc{\cP}{{\mathcal P}}
\nc{\cQ}{{\mathcal Q}}
\nc{\cR}{{\mathcal R}}
\nc{\cS}{{\mathcal S}}
\nc{\cT}{{\mathcal T}}
\nc{\cU}{{\mathcal U}}
\nc{\cV}{{\mathcal V}}
\nc{\cW}{{\mathcal W}}
\nc{\cX}{{\mathcal X}}
\nc{\cY}{{\mathcal Y}}
\nc{\cZ}{{\mathcal Z}}

\nc{\bA}{{\mathbb A}}
\nc{\bB}{{\mathbb B}}
\nc{\bC}{{\mathbb C}}
\nc{\bD}{{\mathbb D}}
\nc{\bE}{{\mathbb E}}
\nc{\bF}{{\mathbb F}}
\nc{\bG}{{\mathbb G}}
\nc{\bH}{{\mathbb H}}
\nc{\bI}{{\mathbb I}}
\nc{\bJ}{{\mathbb J}}
\nc{\bK}{{\mathbb K}}
\nc{\bL}{{\mathbb L}}
\nc{\bM}{{\mathbb M}}
\nc{\bN}{{\mathbb N}}
\nc{\bO}{{\mathbb O}}
\nc{\bP}{{\mathbb P}}
\nc{\bQ}{{\mathbb Q}}
\nc{\bR}{{\mathbb R}}
\nc{\bS}{{\mathbb S}}
\nc{\bT}{{\mathbb T}}
\nc{\bU}{{\mathbb U}}
\nc{\bV}{{\mathbb V}}
\nc{\bW}{{\mathbb W}}
\nc{\bX}{{\mathbb X}}
\nc{\bY}{{\mathbb Y}}
\nc{\bZ}{{\mathbb Z}}

\begin{document}
\title{Stability Manifold of $\bP^{1}$}
\author{So Okada}
\address{Dept of Math\&Stat, UMass Amherst, MA 01003-9305 US, okada@math.umass.edu}
\maketitle
\begin{abstract}
 We describe the stability manifold of the bounded derived category
 $\D(\bP^{1})$ of coherent sheaves on $\bP^{1}$, denoted
 $\Stab(\D(\bP^{1}))$.  This is the first complete picture of a
 stability manifold for a non-Calabi-Yau manifold.
\end{abstract}
\section{Introduction}
T. Bridgeland defined the notion of the stability manifold of a
triangulated category \cite{BRD}, motivated by M. Douglas's work on
$\Pi$-stability of $D$-branes \cite{DOG1}, \cite{DOG2}, and \cite{DOG3}.
Some stability manifolds occur as an approximation of
a part of  moduli of (2,2) superconformal field theories of interest
in algebraic geometry.  Another reason to be interested in the stability
manifold $\Stab(\cT)$ of a triangulated category $\cT$ is that it can be
viewed as a tool for navigation in $\cT$.  For instance, since the
notion of a stability condition refines the notion of a heart,
$\Stab(\cT)$ decomposes into ``cells'' $\Stab_\cC(\cT)$ corresponding to
interesting hearts $\cC$.  Moreover, the extra structure contained in a
stability condition provides  mechanisms of rotation (the action of
$\bC$ in Definition~\ref{def:C action}) and wall crossing (Proposition
\ref{prop:kronecker_wall} and Lemmas
\ref{lem:neighbor_1}-\ref{lem:neighbor_2}), that allow one to
systematically construct new hearts from the known ones.  Throughout
this paper, by a {\it heart} we mean the heart of a bounded
$t$-structure of $\cT$ (the heart actually determines the 
corresponding bounded $t$-structure \cite[Section 3]{BRD}).

 In this paper, we study in some detail the stability manifold for the
 bounded derived category $\D^{b}(\bP^{1})$ of coherent sheaves on
 $\bP^{1}$.  In particular:
 \begin{thm}\label{thm:main}
  $\Stab(\D(\bP^{1}))$ is isomorphic to $\bC^{2}$ as a complex
  manifold.
 \end{thm}
 The strategy is to show that the quotient of $\Stab(\D(\bP^{1}))$ for a
 certain action of $\bC \times \bZ$ is isomorphic to $\bC^{*}$.
 The main technical step is the following list of stability conditions
 of $\D(\bP^{1})$ (the notions used in Theorem \ref{thm:stab} are
 explained in Section \ref{sec:basic-defin-tools}).
  \begin{thm}\label{thm:stab}
   Up to the action of $\Aut(\D(\bP^{1}))$, for any stability condition
   in $\Stab(\D(\bP^{1}))$, there exists some $p>0$ such that
   $\cO(-1)[p]$ and $\cO$ are semistable and $\phi(\cO(-1)[p]),
   \phi(\cO)\in (r, r+1]$ for some $r\in \bR$.

   If $\phi(\cO(-1)[1])<\phi(\cO)$, the multiples of the shifts of
   $\cO(-1)$ and $\cO$ are the only semistable objects. If
   $\phi(\cO(-1)[1])\geq \phi(\cO)$, then all line bundles and torsion
   sheaves are semistable.
  \end{thm}
  
  In Section \ref{sec:basic-defin-tools}, we briefly explain
  parts of \cite{BRD}.  The reader can consult \cite{GM}
  for basic notions of  triangulated categories and  hearts.  In
  Section \ref{sec:stab-cond-dbp1}, we prove Theorem \ref{thm:stab} and
  present all hearts that appear in $\Stab(\D(\bP^{1}))$. In Section
  \ref{sec:stability-manifold}, we find a fundamental domain of
  $\Stab(\D(\bP^{1}))/\bZ\times \bC$ in Lemma \ref{lem:3}, and prove
  Theorem \ref{thm:main}.  In Section \ref{sec:walls-hearts-stabdbp}, we
  explicitly describe how $\Stab(\D(\bP^{1}))$ is glued from ``cells''
  corresponding to interesting hearts of $\D(\bP^{1})$.

\

The author thanks his adviser I. Mirkovi\'c for suggesting the problem
and for countless advice, and thanks R. Kusner, E. Markman, E. Cattani,
and T. Braden for discussions.

\section{Basic definitions and tools}\label{sec:basic-defin-tools}
\subsection{Definition of stability conditions}\label{subsub:stab}
  Let us define the stability conditions using the 
  following two ``filtrations''
  for our model
  (we will loosely use the word {\it filtration} for sequences of
  exact triangles as below).  First, a heart $\cA$ of
  a triangulated category $\cT$ gives a filtration of each object $E\in
  \cT$. For example, the standard heart $\Coh \bP^{1}$ of $\D(\bP^{1})$
  gives the following filtration for each object $E\in \D(\bP^{1})$ by
  taking $E_{k}=\tau_{\leq k} E$.
\[
\xymatrix@C=.5em{
0_{\ }\ar@{=}[r]&E_m\ar[rr]&&E_{m+1}\ar[rr]\ar[dl]&&E_{m+2}\ar[rr]\ar[dl]&&\ldots\ar[rr]&&E_{n-1}\ar[rr]&&E_n\ar[dl]\ar@{=}[r]&E_{\ }\\
&& A_{m+1}\ar@{-->}[ul]&&A_{m+2}\ar@{-->}[ul]&&&&&&A_n \ar@{-->}[ul] 
}
\]
Each cone $A_{j}=H^{j}(E)[-j]$ lies in $\Coh(\bP^{1})[-j]$.

  We refine the filtration above. Now, $A_{n}[n]=I\oplus L$ for a
  torsion sheaf $I$ and $L=\cO(s_{2})^{t_{2}}\oplus \cdots \oplus
  \cO(s_{u+1})^{t_{u+1}}$ for some $s_{2}>\cdots >s_{u+1}$.  Then
\[
\xymatrix@C=.5em{
E_{n-1
}\ar@{=}[r]&E^{0}_{n}\ar[rr]&&E^{1}_{n}\ar[rr]\ar[dl]&&E^{2}_{n}\ar[rr]\ar[dl]&&\ldots\ar[rr]&&E^{u}_{n}\ar[rr]&&E_n\ar[dl]\ar@{=}[r]&
E_{\ }\\
&& A^1_{m+1}\ar@{-->}[ul]&&A^{2}_{m+2}\ar@{-->}[ul]&&&&&&A^{u+1}_n \ar@{-->}[ul] 
}
\]
for $A^{1}_{n}=I[-n]$, $A^{2}_{n}=\cO(s_{2})^{\oplus t_{2}}[-n], \cdots,
A^{u+1}_{n}=\cO(s_{u+1})^{\oplus t_{u+1}}[-n]$.  

The crucial observation is the following; the property that the
rightward $\Hom$ is zero for cones in the first filtration remains true
in the second filtration; i.e., $\Hom_{\cT}(A_{i}, A_{j})=0$ for $i<j$,
and $\Hom_{\Coh \bP^{1}[-n]}(A_{n}^{i}, A_{n}^{j})=0$ for $i<j$.  This
situation is axiomatized in the following definition.
\begin{defin}\cite[Definition 1.1]{BRD}\label{def:stab}
A {\it stability condition} $(Z, \cP)$ on a triangulated category
$\cT$ consists of a group homomorphism $Z: K(\cT)\to\bC$ called the
{\it central charge}, and a family of full additive subcategories
$\cP(\phi)$ of $\cT$ indexed by real numbers $\phi$, called the {\it
slicing}, with the following properties:
\begin{itemize}
\item[(a)] if $E\in \cP(\phi)$ then $Z(E)=m(E)\exp(i\pi\phi)$ for some
	   $m(E)\in\bR_{>0}$;
 \item[(b)] for all $\phi\in\bR$, $\cP(\phi+1)=\cP(\phi)[1]$;
 \item[(c)] if $\phi_1>\phi_2$ and $A_j\in\cP(\phi_j)$ 
	   then $\Hom_{\cT}(A_1, A_2)=0$;
 \item[(d)] for each nonzero object $E\in\cT$ there is a finite 
	   sequence of real numbers
	   \[\phi_1>\phi_2> \cdots >\phi_n\]
	   and a collection of triangles
	   \[
	   \xymatrix@C=.5em{
	   0_{\ }\ar@{=}[r]&E_{0}\ar[rr]&&E_{1}\ar[rr]\ar[dl]&&E_{2}\ar[rr]\ar[dl]&&\ldots\ar[rr]&&E_{n-1}\ar[rr]&&E_n\ar[dl]\ar@{=}[r]&
	   E_{\ }\\
	   && A_1\ar@{-->}[ul]&&A_{2}\ar@{-->}[ul]&&&&&&A_{n} \ar@{-->}[ul] 
	   }
	   \]
	   with $A_j\in\cP(\phi_j)$ for all $j$.
\end{itemize}
 We call the filtration the {\it Harder-Narashimhan filtration} (or {\it
 HN filtration} for short) of $E$ with respect to $(Z, \cP)$, and if
 $E\in \cP(\phi)$ then we call $\phi$ the {\it phase} of $E$.  
 HN-filtrations are unique (when they exist).  
\end{defin}
\subsection{Definition and local structure of $\Stab(\cT)$}
For an interval $I$, let $\cP(I)$ be the full subcategory of $\cT$
generated under extensions by $\cP(\phi)$ for all $\phi\in I$ (we say that $B$
is an extension of $A$ and $C$ in $\cT$ if there exists an exact
triangle $A\to B \to C$ in $\cT$).  A slicing $\cP$ of a triangulated
category $\cT$ is said to be {\it locally-finite}, if for any $t\in \bR$
we have an open interval $I$ around $t$, such that each object in
$\cP(I)$ has a finite-length Jordan-H$\ddot{\text{o}}$lder filtration.
A stability condition $(Z, \cP)$ is locally-finite if the corresponding
slicing $\cP$ is.
\begin{defin}\cite[Section 6]{BRD}
 For a triangulated category $\cT$, $\Stab(\cT)$ is the set of all
 locally-finite stability conditions on $\cT$.
T. Bridgeland defines a topology on $\Stab(\cT)$.
Its characterizing property is
\end{defin}
\begin{thm}\label{thm:local}\cite[Theorem 1.2]{BRD}
 Let $\cT$ be a triangulated category.  For each connected component
 $\Sigma\subset\Stab(\cT)$ there is a linear subspace $V(\Sigma)\subset
 (\cK(\cT)\otimes\bC)^*$ with a well-defined linear topology, and such
 that the map $\cZ:\Sigma\to V(\Sigma)$, which maps a stability condition
 $(Z, \cP)$ to its central charge $Z\in V(\Sigma)$, is a local homeomorphism.
\end{thm}
Hence, when $K(\cT)$ has finite rank, $\Stab(\cT)$ is a complex
manifold, called the {\it stability manifold} of $\cT$.
\subsection{Hearts and stability conditions}\label{subsub:hearts}
Let us see the relation between hearts and stability conditions.  By a
heart of $\cT$, we mean a heart of a bounded $t$-structure on $\cT$. For
a given stability condition $(Z, \cP)$, $\cP((r, r+1])$ is a heart for
any $r\in \bR$, since a heart can be characterized as follows.
\begin{lemma}\label{lem:heart}\cite[Lemma 3.2]{BRD}
 Let $\cA\subset\cT$ be a full additive subcategory of a triangulated
 category $\cT$.  Then $\cA$ is a heart of $\cT$ if and only if the
 following two conditions hold:
\begin{itemize}
 \item[(a)]if $k_1>k_2$ are integers and $A, B\in \cA$ 
	   then $\Hom_{\cT}(A[k_1], B[k_2])=0$;
 \item[(b)]for every nonzero object $E\in\cT$ there is a finite
	   sequence of integers
	   \[k_1>k_2>\cdots>k_n\]
	   and a collection of triangles
	   \[
	   \xymatrix@C=.5em{
	   0_{\ }\ar@{=}[r]&E_{0}\ar[rr]&&E_{1}\ar[rr]\ar[dl]&&E_{2}\ar[rr]\ar[dl]&&\ldots\ar[rr]&&E_{n-1}\ar[rr]&&E_n\ar[dl]\ar@{=}[r]&
	   E_{\ }\\
	   && A_1\ar@{-->}[ul]&&A_{2}\ar@{-->}[ul]&&&&&&A_{n} \ar@{-->}[ul] 
	   }
	   \]
	   with $A_j\in\cA[k_j]$ for all $j$.
\end{itemize}
\end{lemma}
The main relation between hearts and stability conditions is the
following.
 \begin{prop}\label{prop:extend}\cite[Proposition 5.3]{BRD}
  To give a stability condition on a triangulated category $\cT$ is
  equivalent to giving a heart on $\cT$ and a centered slope-function on
  it with the Harder-Narasimhan property.
 \end{prop}
 Here, a {\it centered slope-function} $Z$ on a heart $\cA$ is a group
 homomorphism $Z:K(\cA)\to \bC$, such that for $0 \neq E \in \cA$,
 $Z(E)$ lies in $H\stackrel{\mbox{{\tiny
 def}}}{=}\{r\exp(i\pi\phi):r>0\text{ and }0<\phi\leq 1\}\subset \bC$.
 We define the {\it slope} of $E\neq 0$, denoted by $\phi(E)$, to be
 \begin{align*}
  \phi(E)=\frac{1}{\pi} \arg Z(E)\in (0, 1].  
 \end{align*} 
 Object $0\neq E\in\cA$ is called {\it semistable} if for any subobject
 $0\neq A$ of $E$  $\phi(A)\leq\phi(E)$.  Moreover, $Z$ is said
 to have the {\it Harder-Narasimhan property} (or {\it HN-property} for
 short), if for every nonzero object $E$ of $\cA$ there is a finite
 short exact sequences (which we draw as triangles) in $\cA$
	   \[
	   \xymatrix@C=.5em{
	   0_{\ }\ar@{=}[r]&E_{0}\ar[rr]&&E_{1}\ar[rr]\ar[dl]&&E_{2}\ar[rr]\ar[dl]&&\ldots\ar[rr]&&E_{n-1}\ar[rr]&&E_n\ar[dl]\ar@{=}[r]&
	   E_{\ }\\
	   && A_1\ar@{-->}[ul]&&A_{2}\ar@{-->}[ul]&&&&&&A_{n} \ar@{-->}[ul] 
	   }
	   \]
with $\phi(A_{i})>\phi(A_{j})$ for $i<j$.

Proposition \ref{prop:extend} says that a centered slope-function on a
heart with the HN-property can be extended to a stability condition.
Notice that for a heart $A$ of $\cT$, we always have $K(A)=K(\cT)$. So,
such an extension is achieved just by setting the central charge to be the
same to the centered slope-function, and by setting the slice
$\cP(\psi+k)$, for $\psi\in (0, 1]$ and $k\in \bZ$, to be the full
additive subcategory of $\cT$ consisting of $E[k]$ with $\phi(E)=\psi$.
A non-zero object $E\in \cP(\phi)$ for each $\phi\in \bR$ is also
called semistable.

Let us see how it works in our example $\Coh \bP^{1}$.  
Notice that we have to require
 $Z(\cO_{x})\in \bR_{<0}$ and $Z(\cO)\in H \setminus \bR_{<0}$ 
in order to have
 $Z(\Coh\bP^{1})\subset H$.  Let us graphically present a centered
 slope-function $Z$ on the standard heart $\Coh \bP^{1}$ by the
 following figure.
  \begin{figure}[H]
   \begin{center}
    \scalebox{.7}{
    \input{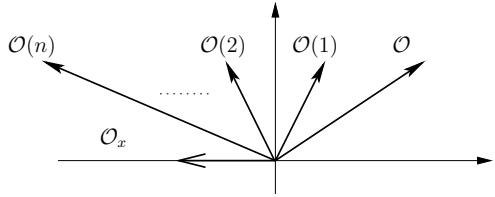}
    }
    \caption{A centered slope-function $Z$ on $\Coh \bP^{1}$}
    \label{fig:std}
   \end{center}
  \end{figure}  
  Then the second filtration in \ref{subsub:stab} proves that
  $Z$ has the HN-property, with all the line bundles and the torsion
  sheaves semistable.  Hence $Z$ extends to a stability condition $(Z,
  \cP)$.

  \subsection{Wall crossing and rotation}
  Let $\mathrm{K}$ be the {\it Kronecker quiver}
  $(\cdot\rightrightarrows\cdot)$, and $\Rep(\mathrm{K})$ be the
  category of the representations of the quiver.  We will construct
  Beilinson's {\it Kronecker heart} $\cB$ in $\D(\bP^{1})$ (see
  \cite{BEI}), and some stability conditions with the heart $\cB$.  We
  will use two operations on $\Stab(\cT)$; {\it rotation} affects the
  heart but preserves the semistable objects, while {\it wall crossing}
  fixes the heart and affects semistable objects.

    For $S\subset \cT$, we let $\langle S \rangle$ be the full subcategory
    generated under extension by  objects in $S$.  {\it Rotation} and {\it
    rescaling} of a stability condition are the imaginary and the real
    parts of the $\bC$-action below (it will be discussed in Proposition
    \ref{prop:rotation}).
    \begin{defin}\label{def:C action}
     Let $(Z, \cP)$ be a stability condition and $z=x+i y\in \bC$.  Then
     $z*(Z, \cP)$ is defined to be $z*Z=e^{z}Z$ and
     $(z*\cP)(\phi)=\cP(\phi-y/\pi)$.
    \end{defin}
   
   \begin{prop}\label{prop:kronecker}
    Let $(Z, \cP)$ be a stability condition as in the end of the Section
    \ref{subsub:hearts}.  Then rotation by $z=-i \phi(\cO)$ gives a
    stability condition $z*(Z, \cP)=(\overline{Z}, 
    \overline{\cP})$ such that the
    heart $\overline{\cP}((0,1])$ is equivalent to $\Rep(\mathrm{K})$ and all
    line bundles and torsion sheaves are semistable.
   \end{prop}
   \begin{proof}
    We easily check that $\overline{\cP}((0,1]) =\langle \cO, \cO(-1)[1]
    \rangle$, since $\overline{\cP}(\cO(-1))=0$.  Let $P=\cO\oplus \cO(1)$.
    The functor $\mathrm{R}\Hom(P, -)$ sends $\cO$ and
    $\cO(-1)[1]$ to the irreducible generators of $\Rep(\mathrm{K})$,
    and gives the equivalence between the categories $\overline{\cP}((0,1])$
    and $\Rep(\mathrm{K})$ (see \cite[Section 6]{BON}).
   \end{proof} 

   Now we fix the heart $\cB$ and vary the centered slope-function
   $\widetilde{Z}$ on it.  Since $\cB$ has finite-length objects, centered
   slope-functions with HN-property are defined by any choice of the
   values $\widetilde{Z}(\cO), \widetilde{Z}(\cO(-1)[1])\in H$, on irreducible
   objects $\cO, \cO(-1)[1]$ in $\cB$.
  
   \begin{prop}\label{prop:kronecker_wall}
    If $\widetilde{\phi}(\cO)<\widetilde{\phi}(\cO(-1)[1])$, then we get
    $\overline{Z}$ above.  
    If $\widetilde{\phi}(\cO)=\widetilde{\phi}(\cO(-1)[1])$,
    then any non-zero object in $\widetilde{\cP}((0,1])$ is semistable.  If
    $\widetilde{\phi}(\cO)>\widetilde{\phi}(\cO(-1)[1])$, then the multiples of
    $\cO$ and $\cO(-1)[1]$ are the only semistable objects in
    $\widetilde{\cP}((0,1])$.
   \end{prop}
   \begin{proof}
    If $\widetilde{\phi}(\cO)=\widetilde{\phi}(\cO(-1)[1])$, then all non-zero
    objects in $\cB=\langle \cO, \cO(-1)[1] \rangle$ lie in the same
    slope.

    Consider the last case.  We see that $\cO(n)$ is not semistable for
    $n>0$, since we have the triangle $\cO\to \cO(n) \to I$ for some
    torsion sheaf $I$ and $\widetilde{\phi}(\cO)>\widetilde{\phi}(\cO(n))$.
    Similarly, $\cO(-n-1)[1]$ for $n>0$ is not semistable because of the
    triangle $I \to \cO(-n-1)[1]\to \cO(-1)[1]$ for some torsion
    sheaf $I$, and any torsion sheaf is not semistable because of
    the triangle $\cO \to \cO_{x} \to \cO(-1)[1]$ for any $x$.
   \end{proof}

   T. Bridgeland defined {\it wall}, as a codimension-one submanifold of
   a stability manifold such that as one varies a stability condition, a
   semistable object can only become non-semistable if one crosses a wall
   (see \cite[Section 8]{BRD1}).  From Proposition
   \ref{prop:kronecker_wall}, $W=\{ (Z, \cP)\in \Stab(\D(\bP^{1}))\mid
   \phi(\cO)=\phi(\cO(-1)[1])\}$ is a wall.
   
 \section{Stability conditions for  $\D(\bP^{1})$}\label{sec:stab-cond-dbp1}
 The remark after \cite[Theorem 3.1]{BON1} says that $\Aut (\D(\bP^{1}))\cong
 \Aut \bP^{1} \ltimes (\mathrm{Pic} (\bP^{1}) \oplus \bZ)$, where $\bZ$ is
 generated by the shift $[ 1 ]$.  We frequently use
 $\mathrm{Pic}(\bP^{1})\oplus \bZ\cong \bZ \oplus \bZ$, while $\Aut
 \bP^{1}$ acts trivially on $\Stab(\D(\bP^{1}))$.

  One way to obtain Theorem \ref{thm:stab} is to classify the centered
  slope-functions with HN-property for each heart in the list provided
  by \cite[Theorem 6.12]{GKR}.  However, for our purpose it is simpler
  to use the following lemma.
  \begin{lem}\cite[Lemma 6.6]{GKR}\label{lem:GKR}
   \begin{enumerate}
    \item[(a)] For $n\in \bZ$, there exist exact triangles
	       \begin{align}
		&\cO(k+1)^{\oplus n-k}		
		&\to&
		\cO(n)
		\to
		\cO(k)^{\oplus n-k-1}[1]&		
		\ \ &\mbox{if } n>k+1,\label{eq:a-1}\\ 
		&\cO(k+1)^{\oplus k-n}[-1]
		&\to&
		\cO(n)
		\to
		\cO(k)^{\oplus k-n+1}&
		\ \ &\mbox{if } n<k.\label{eq:a-2}
	       \end{align}
    \item[(b)] For $x\in \bP^{1}$ and $k\in \bZ$, 
	       there exists an exact triangle
	       \begin{align}
		\cO(k+1)
		&\to
		\cO_{x}
		\to
		\cO(k)[1]. \ \ \label{eq:b-1}
		\end{align}
    \item[(c)] Any triangle $A\to M \to B$ with $\Ext^{\leq 0}(A, B)=0$
	       and $M$ either $\cO(n)$ or $\cO_{x}$ is in the form of
	       $(a)$ or $(b)$.
    \item[(d)]\label{lem:key-d} If some line bundle or torsion sheaf is
	       not semistable, then there exist $k, n\in
	       \bZ$ such that the shifts of $\cO(k)$ and $\cO(k+1)$ are
	       semistable and the triangle found in $(a)$ or $(b)$ is
	       the HN-filtration.
   \end{enumerate}
  \end{lem}
 \begin{proof}
  For (a)-(c), see the proof after \cite[Remark 6.8]{GKR}.  For (d),
  see the proof of the semistability of $\cO(k)$ and
  $\cO(k+1)$ after the statement of \cite[Lemma 6.6]{GKR}.
  \end{proof}

  \begin{cor}\label{cor:phases}
   If there exists $k$ such that $\cO(k)$ and $\cO(k+1)$ are semistable
   and $\phi(\cO(k+1))> \phi(\cO(k)[1])$, then no line bundle or torsion
   sheaf is semistable except $\cO(k)$ and $\cO(k+1)$.
  \end{cor}
  \begin{proof}
   If $\cO(n)$ for $n>k+1$ is semistable, then $\cO(n)$ itself is the
   HN-filtration, but the triangle \eqref{eq:a-1} is also the
   HN-filtration, since $\cO(k+1)$ and $\cO(k)[1]$ are semistable and
   the cone phases are decreasing by the assumption.  So, by the
   uniqueness of the HN-filtrations, $\cO(n)$ cannot be semistable.
   Likewise, $\cO(n)$ for $n<k$ and torsion sheaf are not semistable,
   because of the triangle \eqref{eq:a-2} and the triangle
   \eqref{eq:b-1} respectively.
  \end{proof}

 \def\proofname{Proof of Theorem \ref{thm:stab}}
 \begin{proof}
  Consider the case when there exists some non-semistable line bundle.
  Lemma \ref{lem:GKR} (d) says that there is a line bundle $\cO(n)$ with
  the HN-filtration
  \[\xymatrix@C=.5em{\cO(k+1)^{\oplus p}[j] \ar[rr] && \cO(n)  \ar[dl] \\
  & \cO(k)^{\oplus s}[j+1] \ar@{-->}[ul]} \]
  for some $k \in \bZ$, $j\in \{0, 1\}$, $p, s>0$. 
  The decreasing property of the phases the cones of the HN-filtration
  implies
  \begin{align}
   \phi(\cO(k+1))> \phi(\cO(k)[1]). \label{ineq:case-1}
  \end{align}
  Hence, Corollary \ref{cor:phases} says that up to shifts $\cO(k)$ and
  $\cO(k+1)$ are the only semistable objects in $\D(\bP^{1})$; every
  object in $\D(\bP^{1})$ can be decomposed into a direct sum of some
  shifts of line bundles and torsion sheaves, since the homological
  dimension of $\Coh \bP^{1}$ is one.
  After tensoring with $\cO(-k-1)$, the inequality \eqref{ineq:case-1}
  gives $\phi(\cO)> \phi(\cO(-1)[1])$ and hence
  \begin{align*}
   \phi(\cO(-1)[p]), \phi(\cO)&\in (r, r+1]   
  \end{align*}
  for some $r\in \bR$ and $p>0$.

  Consider the case when all line bundles are semistable.  Since there
  are more than two line bundles that are semistable, Corollary
  \ref{cor:phases} says that all torsion sheaves have to be semistable,
  and in particular
  \begin{align}
   \phi(\cO(k+1))\leq \phi(\cO(k)[1])   \label{ineq:case-2-2}
  \end{align}  
  for any $k$.

  Let us see that there exists $k$ such that
  \begin{align}
   \phi(\cO(k)[1])-1<\phi(\cO(k+1)). \label{ineq:case-2-3}
  \end{align}
  If not, $\phi(\cO(k)[1])-1\geq \phi(\cO(k+1))$ for all $k$; i.e.,
  $\phi(\cO(k)[1])\geq \phi(\cO(k+1)[1])$ for all $k$.  Since $\Hom(
  \cO(k)[1], \cO(k+1)[1])\neq 0$, $\phi(\cO(k)[1])=\phi(\cO(k+1)[1])$
  for all $k$ by Property (c) in Definition \ref{def:stab}.  However,
  the phases of the line bundles cannot be the same, since
  $Z(\cO(k+1))-Z(\cO(k))=Z(\cO_{x})\neq 0$ for all $k$ and they cannot
  be on the same ray.  Therefore, the inequality \eqref{ineq:case-2-3}
  holds for some $k$, and
  \begin{align*}
   0\leq \phi(\cO(-1)[1])-\phi(\cO)<1
  \end{align*}
  up to $\Aut(\D(\bP^{1}))$.
 \end{proof}  
  \def\proofname{Proof}
  
  Theorem \ref{thm:stab} only tells us necessary conditions for the
  stability conditions in $\Stab(\D(\bP^{1}))$, but Proposition
  \ref{prop:exp} says all cases in the theorem actually exist and
  reveals all hearts on which we can define a centered slope-function for
  each case.

  \begin{prop}\label{prop:exp}
   For any $\alpha, \beta\in \bR$ such that $\alpha>\beta-1$, and any
   $m_{\alpha}, m_{\beta} \in \bR_{>0}$, there exists a unique stability
   condition $(Z, \cP)$ such that $\phi(\cO(-1)[1])=\beta$ and
   $\phi(\cO)=\alpha$, and $Z(\cO)=m_{\alpha}e^{i \pi \alpha}$ and
   $Z(\cO(-1)[1])=m_{\beta}e^{i \pi \beta}$.

   Furthermore, we have the following cases:
   \begin{enumerate}
    \item if $\alpha>\beta$, then for any $r\in \bR$, there exist
	  $p,q\in \bZ$ such that $\cP((r-1, r])=\langle\cO(-1)[p+1],
	  \cO[q]\rangle$ and $p-q\in (\alpha-\beta-1,\alpha-\beta+1)$;
    \item if $\alpha< \beta$, then for any $r$, either there exist
	  $i,j\in \bZ$ such that $\cP((r-1, r])=\langle \cO(i-1)[1+j],
	  \cO(i)[j]\rangle$ and $\phi(\cO(i-2)[1+j])> r \geq
	  \phi(\cO(i-1)[1+j])$, or there exists $j\in \bZ$ such that
	  $\cP((r-1,r])=\Coh \bP^{1}[j]$ and $r=\phi(\cO_{x}[j])$;
    \item if $\alpha=\beta$, then for any $r$, there exists $j$ such
	  that $\langle \cO(-1)[1+j], \cO[j]\rangle =\cP((r-1, r])$.
   \end{enumerate}
  \end{prop}
  \begin{proof}
   In all cases, for any $k\in \bZ$, $\cP$ has to satisfy
   $\cP(\alpha+k)=\langle \cO[k] \rangle$, $\cP(\beta+k)=\langle
   \O(-1)[1+k] \rangle[k]$. 
   
   Consider the first case.  By Theorem \ref{thm:stab}, the
   multiples of the shifts of $\cO$ and $\cO(-1)[1]$ are
   the only semistable objects.
   So we have to put $\cP(\phi)=0$ for all $\phi\in \bR$ except when
   $\phi=\alpha+k$ or $\beta+k$ for some $k\in \bZ$.  Then $(Z, \cP)\in
   \Stab(\D(\bP^{1}))$; (c) in Lemma \ref{lem:GKR} implies (c) in
   Definition \ref{def:stab}: the triangles \eqref{eq:a-1},
   \eqref{eq:a-2}, and \eqref{eq:b-1} are the HN-filtrations for all
   summands of each object in $\D(\bP^{1})$.  If $\alpha-\beta\in \bZ$,
   for any $r\in \bR$, there exists $p, q\in \bZ$ such that $\cP((r-1,
   r]) =\langle \cO(-1)[p+1], \cO[q] \rangle$ and $p-q=\alpha-\beta$. If
   $n\not\in \bZ$, then there exists $p, q\in \bZ$ $\cP((r-1,r])=\langle
   \cO(-1)[1+p], \cO[q] \rangle$ and $p-q$ is either the smallest or the
   largest possible integer to $\alpha-\beta$.
   
   Consider the next case.  By Theorem \ref{thm:stab}, all line bundles
   and torsion sheaves are semistable. So for all $\phi\in (0,1]$ and
   $k\in \bZ$, we have to put $\cP(\phi+k)=\langle \cO(n)[k] \rangle$
   when $Z(\cO(n))=m(\cO(n))e^{i\pi \phi}$ for some $n>0$ and
   $m(\cO(n))\in \bR_{>0}$, $\cP(\phi+k)=\langle \cO(-n-1)[1+k] \rangle$
   when $Z(\cO(-n-1))=m(\cO(-n-1))e^{i\pi \phi}$ for some $n>0$ and
   $m(\cO(-n-1))\in \bR_{>0}$, $\cP(\phi+k)=\langle \cO_{x}[k]\mid x\in
   \bP^{1} \rangle$ when $Z(\cO_{x})=m(\cO_{x})e^{i\pi \phi}$ for some
   $x\in \bP^{1}$ and $m(\cO_{x})\in \bR_{>0}$, and $\cP(\phi+k)=0$ for
   the other cases.  Notice that $Z(\cO_{x})=m(\cO_{x}) e^{i \pi \psi}$
   for some $\beta>\psi>\alpha$ and $m(\cO_{x})\in \bR_{>0}$ because of
   the triangle $\cO\to \cO_{x} \to \cO(-1)[1]$.  So for
   $z=e^{i\pi(1-\phi)}$, $z*Z$ gives us a centered slope-function with
   HN-property in \ref{subsub:hearts}. Hence we have $z*(Z, \cP)\in
   \Stab(\D(\bP^{1}))$, and $(Z, \cP)\in \Stab(\D(\bP^{1}))$.

   Finally, let $\cP(\alpha)=\langle \cO(-1)[1], \cO \rangle$.  Then
   $(Z, \cP)\in\Stab(\D(\bP^{1}))$ by Proposition
   \ref{prop:kronecker_wall}.
  \end{proof}

  \begin{cor}\label{cor:all_hearts}
   The hearts on which we can impose a centered slope-function with
   HN-property are $\cC_{j}\stackrel{{\tiny {\rm def}}}{=}\Coh
   \bP^{1}[j]$ and $\cC_{p, i, j}\stackrel{{\tiny {\rm def}}}{=}\langle
   \cO(i-1)[p+j], \cO(i)[j] \rangle$ for all $i, j\in \bZ$ and $p> 0$.
  \end{cor}
  \begin{proof}
   The hearts found in Proposition \ref{prop:exp} are
   $\langle\cO(i-1)[p], \cO(i)[q]\rangle$ for $i\in \bZ$ and any $p,q\in
   \bZ$ such that $p-q>0$, and $\Coh \bP^{1}[j]$ for any $j$.  Any heart
   on which we can put a centered slope-function with HN-property lives
   in their orbits under the action of $\Aut(\D(\bP^{1}))$ by Theorem
   \ref{thm:stab}.
  \end{proof}

  \begin{rmk}\label{rmk:grd}
   The subcategory $\cA=\langle \cO_{x}, x\in \bP^{1}; \cO(n)[1], n\in
   \bZ\rangle$ is a heart because it is the image of $\Coh \bP^{1}$
   under the Grothendieck duality functor $\bD(-)={\rm
   R}\mathcal{H}om(-, \omega_{\bP^{1}})[1]$.  However, it carries no
   centered slope-function with HN-property by Corollary
   \ref{cor:all_hearts}.
   This shows that the notion of stability is not invariant under
   passing from $\cT$ to $\cT^{\rm opp}$.  So if we call the ones in
   Definition~\ref{def:stab} {\it right stability conditions}, there is
   a notion of {\it left stability conditions} on $\cT^{\rm opp}$, and
   we have such left stability conditions on $\cA\subset
   \D(\bP^{1})^{\rm opp}$.

   For any $P\subset \bP^{1}$, a subcategory $\cA(P)=\langle \cO_{x},
   x\in P; \cO_{y}[1], y\not\in P; \cO(n)[1], n\in \bZ\rangle$ is a
   heart (notice that $\cA(P)=\cA$).  Up to $\Aut(\D(\bP^{1}))$, all
   hearts that do not bear any centered slope-function with HN-property
   are $\cA(P)$ by Corollary \ref{cor:all_hearts} and \cite[Theorem
   6.12]{GKR}.
  \end{rmk}
  
 \section{Stability Manifold}\label{sec:stability-manifold}
 \subsection{Quotient $\Stab(\cT)/\bC$}
  T. Bridgeland observes that $\widetilde{\GL^{+}}(2, \bR)$ (the
  universal covering of $\GL^{+}(2, \bR)$) acts on $\Stab(\cT)$ for any
  triangulated category $\cT$ (\cite[Lemma 8.2]{BRD}).  We notice that
  the $\bC$-action in Definition \ref{def:C action} is a holomorphic
  part of this action.
  
 \begin{prop}\label{prop:rotation}
  The $\bC$-action is holomorphic, free, coincides with the action of a
  subgroup of $\widetilde{\GL^{+}}(2, \bR)$, and contains the shifts.
  The quotient $\Stab(\cT)/\bC$ is a complex manifold,
  modeled on a projective space of a topological vector space.
 \end{prop}
 \begin{proof}
  By Theorem \ref{thm:local}, holomorphicity follows from the
  holomorphicity of the $\bC$-action on a vector space via
  multiplication by $e^{z}$.  The stabilizers are trivial since $z *
  Z=Z$ implies $e^{x}=1$, so $x=0$, while $z*\cP=\cP$ gives $y=0$.
  Notice that the $\bC$-action by $z$ has the same effect as the action
  of $(e^{-x}A, y)\in \widetilde{\GL^{+}}(2, \bR)$ on a stability condition
  (notation from \cite[Lemma 8.2]{BRD}), where $A$ is the rotation by
  the angle $-\pi y$. The shift $[1]$ can be realized as the action of
  $i \pi \in \bC$.  The action of $\bC$ on the manifold $\Stab(\cT)$ is
  point-wise free, and
  locally isomorphic to the $\bC^{*}$-action
  on $V(\Sigma)\setminus \{0\}$ (notation from
  Theorem \ref{thm:local}).
  This implies that
  $\Stab(\cT)/\bC$ is a
  manifold as claimed.
 \end{proof}

 \subsection{Quotient $\Stab(\D(\bP^{1}))/(\bZ)\bC$}
  We will denote by $(\bZ)$ the copy of $\bZ$ that acts on
  $\bD(\bP^{1})$ by the tensoring with line bundles.

  In Lemma \ref{lem:1}, we show that Theorem \ref{thm:stab} gives a
  domain $X$ that contains a fundamental domain of
  $\Stab(\D(\bP^{1}))/(\bZ)\bC$,
 but we still need Lemma \ref{lem:3} to shrink 
 $X$ so that we 
 avoid  repetitions with respect to the action of $(\bZ)\bC$.

 \begin{lem}\label{lem:1}
  Let $X$ be the subset of $\Stab(\D(\bP^{1}))$ consisting of all
  stability conditions $(Z, \cP)$ with the following properties: (a)
  $\cO(-1)[1], \cO$ are semistable; (b) $\phi(\cO(-1)[1])=1$ and
  $m(\cO(-1)[1])=1$; (c) $\phi(\cO)>0$.  Then $(\bZ)\bC\cdot
  X=\Stab(\D(\bP^{1}))$, and $X$ is isomorphic to the open upper
  half-plane $\bH$, an isomorphism is given by $\log(m(\cO))+i\pi
  \phi(\cO):X\cong \bH$.
 \end{lem}
 \begin{proof}
  Up to $\Aut(\D(\bP^{1}))$, Theorem \ref{thm:stab} says that, for each
  stability condition $(Z, \cP)$, we have $r\in \bR$ such that $\cO(-1),
  \phi(\cO)$ are semistable and the slope are in $(r, r+1]$.  By the
  action of $\bC$, we can assume $\phi(\cO(-1)[1])=1$ and
  $m(\cO(-1)[1])=1$.  So $(r, r+1]\ni \phi(\cO(-1)[p])=p\geq 1$ forces
  $r>0$, hence $\phi(\cO)>0$.
  
  The slope and the length of $\cO(-1)[1]$ are fixed by (b).  So each
  $(Z, \cP)\in X$ in the stability manifold can be uniquely represented
  by $Z(\cO)$ on the $n$-th sheet of the Riemann surface of $\Log z$,
  where $n$ is the greatest integer such that $\phi(\cO)/2 \geq n$. Here
  $m(\cO), \phi(\cO)>0$ by (a) and (c).
 \end{proof}

 \begin{lem}\label{lem:3}
  A fundamental domain of $\Stab(\D(\bP^{1}))/(\bZ)\bC$ is isomorphic to
  $K\stackrel{{\tiny \mathrm{def}}}{=}\{x+iy\in \bC\mid y>0, \cos y \geq
  e^{-\abs{x}}\}$ as in the shaded domain in the figure below.  When
  passing to $\Stab(\D(\bP^{1}))/(\bZ)\bC$ one identifies points on the
  boundary that have the same imaginary part.  
  {\tiny
  \begin{figure}[H]
   \begin{center}
    \scalebox{.7}{
    \input{K_7.pstex_t}
    }
    \caption{A fundamental domain of $\Stab(\D(\bP^{1}))/(\bZ)\bC$}
    \label{fig:K}
   \end{center}
  \end{figure}  
  }
 \end{lem}
 \begin{proof}
  Notice that the actions by line bundles and by $\bC$ commute.  Let
  $(Z, \cP)\in X$.  First, we will see that if $\phi(\cO)>1$, there are
  no repetitions; i.e., $\bC(\bZ)\cdot (Z, \cP)\cap X=\{(Z, \cP)\}$.  By
  Theorem \ref{thm:stab}, all indecomposable semistable objects are
  shifts of $\cO(-1)$ and $\cO$.  The action of $ \cO(i)\cdot (x+i \pi
  \psi) \in (\bZ) \bC$ changes $(Z, \cP)$ into $(Z', \cP')$ such that
  $\cP'(\phi(\cO)+\psi)=\langle \cO(i) \rangle$ and
  $\cP'(\phi(\cO(-1)[1])+\psi)=\langle \cO(-1+i)[1] \rangle$
  $m'(\cO(i))/m(\cO)=e^{x}$ and $m'(\cO(-1+i)[1])/m(\cO(-1))=e^{x}$.
  Hence, $\cO$ and $\cO(-1)[1]$ are not semistable unless $i=0$, and
  even if $i=0$ we have $\phi'(\cO(-1)[1])\neq 1$ or $m'(\cO(-1)[1])\neq
  1$ unless $\phi=x=0$.

  In the remaining case $0<\phi(\cO)\leq \phi(\cO(-1)[1])=1$;
  repetitions $\bC(\bZ)\cdot (\bZ, \cP)\cap X$ are indexed by
  $\bZ\ni \i\mapsto (Z_{i}, \cP_{i}) =(z_{i}, \cO(i))\cdot (Z, \cP)$.
  Here, $z_{i}=\frac{1}{m(\cO(i-1)[1])}+i \pi(1-\phi(\cO(i-1)[1]))$.
  Let us denote $(\dot{Z}_{i}, \dot{\cP}_{i})=\cO(i)\cdot (Z,
  \cP)$; i.e.,
  \begin{alignat*}{3}
   \dot{Z_{i}}(\cO(-1)[1])&=Z(\cO(i-1)[1]),   
   \ \  &\dot{\phi_{i}}(\cO(-1)[1])&=\phi(\cO(i-1)[1]),\\
   \dot{Z_{i}}(\cO)&=Z(\cO(i)),
   &\dot{\phi_{i}}(\cO)&=\phi(\cO(i)).    
  \end{alignat*}
  We have $z_{i} * (\dot{Z}_{i}, \dot{\cP}_{i})=(Z_{i}, \cP_{i})\in X$, since
  $1-(z_{i}*\dot{\phi}_{i})(\cO)= \phi(\cO(i-1)[1])-\phi(\cO(i))<1$ implies
  $(z_{i}*\dot{\phi}_{i})(\cO)\in (0, 1]$.
  Graphically we can explain actions above as follows.
  \newlength{\minitwocolumn}
  \setlength{\minitwocolumn}{0.45\textwidth}
  \addtolength{\minitwocolumn}{0.3\columnsep}
  \begin{figure}[H]
   \begin{tabular}{c c}
    \begin{minipage}{\minitwocolumn}
     \begin{center}
      \scalebox{.7}{\input{K_1.pstex_t}}
      \caption{$Z$}
      \label{fig:central_charge_Z}
     \end{center}
    \end{minipage}
    &
    \begin{minipage}{\minitwocolumn}
     \begin{center}
       \scalebox{.7}{\input{K_4.pstex_t}}
       \caption{$\dot{Z}_{1}$}
     \end{center}
    \end{minipage}
   \end{tabular}
  \end{figure}
  \begin{figure}[H]
   \begin{center}
   \scalebox{.7}{ 
    \input{K_5.pstex_t}
    }
    \caption{$Z_{1}$}
    \label{fig:central_charge_Z_1}
   \end{center}
  \end{figure} 

  Hence, the question is how to pick one $(Z_{i}, \cP_{i})$.  The first
  step is to require $(\dot{Z}_{i}, \dot{\cP}_{i})$ to have the least
  possible $\dot{\phi}_{i}(\cO(-1)[1])-\dot{\phi_{i}}(\cO)$; it is easy
  to see from simple plane geometry that the minimality is achieved
  exactly when $\dot{Z}_{i}(\cO(-1)[1])$ and $\dot{Z}_{i}(\cO)$ are in
  the strip $S_{Z}$ bounded by the two lines that are perpendicular to
  $Z(\cO_{x})$ and contain the initial or end point of $Z(\cO_{x})$.
  For example in Figures
  \ref{fig:central_charge_Z}--\ref{fig:central_charge_Z_1} $S_{Z}$ is
  indicated by the dotted lines.  Let us show the existence of such
  $(\dot{Z}_{i}, \dot{\cP}_{i})$.  Notice that
  $\dot{Z}_{i}(\cO(-1)[1])\in S_{Z}$ if and only if $\dot{Z}_{i}(\cO)\in
  S_{Z}$, since
  $\dot{Z}_{i}(\cO(-1)[1])+\dot{Z}_{i}(\cO)=\dot{Z}_{i}(\cO_{x})=Z(\cO_{x})$.
  Moreover, we have some $i$ such that $\dot{Z}_{i}(\cO(-1)[1])\in
  S_{Z}$, since
  $\dot{Z}_{i}(\cO(-1)[1])-\dot{Z}_{i-1}(\cO(-1)[1])=Z(\cO_{x})$.

  Graphically, the minimality is achieved precisely when $Z(\cO)$ is in the
  shaded region; $Z(\cO(-1)[1]), Z(\cO)\in S_{Z}$ if and only if
  $Z(\cO)=x+iy$ with $x\leq 1$ and $(x-1/2)^{2}+y^{2}\geq 1/4$.
  \begin{figure}[H]
   \begin{center}
    \scalebox{.7}{
    \input{K_3.pstex_t}
    }
    \caption{Domain  for $Z(\cO)$}
    \label{fig:subdomain_in_H}
   \end{center}
   \end{figure}  

  The remaining repetitions occur precisely when one of $Z(\cO(-1)[1])$
  and $Z(\cO)$ makes a right angle with $Z(\cO_{x})$; i.e., $(Z_{-1},
  \cP_{-1})\in X$ and $Z_{-1}(\cO(-1)[1])\in S_{Z}$, or $(Z_{1},
  \cP_{1})\in X$ and $Z_{1}(\cO(-1)[1])\in S_{Z}$. So two boundary
  points are identified when a ray from the origin connects them; when
  we have a right triangle between $Z(\cO_{x})$ and $Z(\cO)$ as
  $Z(\cO)=1+iy$ for some $y>0$, then
  $Z_{1}(\cO)=\frac{1}{1+y^{2}}(1+iy)$, since $\dot{Z}_{1}(\cO(-1)[1])
  =-1+iy$ and $\dot{Z}_{1}(\cO)=1$.  This situation can be presented graphically
  in the following figures.  
  \setlength{\minitwocolumn}{0.45\textwidth}
  \addtolength{\minitwocolumn}{0.3\columnsep}
  \begin{figure}[H]
   \begin{tabular}{c c}
    \begin{minipage}{\minitwocolumn}
     \begin{center}
      \scalebox{.7}{\input{K_2.pstex_t}}
     \caption{$Z$}
     \end{center}
    \end{minipage}
    &
    \begin{minipage}{\minitwocolumn}
     \begin{center}
     \scalebox{.7}{ 
     \input{K_22.pstex_t}
     }
     \caption{$Z_{1}$}
     \end{center}
    \end{minipage}
   \end{tabular}
  \end{figure}

  Let $X'\subset X$ consist of all $(Z, \cP)\in X$ such that
  $Z(\cO)\in S_{Z}$. Then $\log(m(\cO))+i \pi \phi(\cO)$ gives us an
  isomorphism between $X'$ and $K$; we identify two boundary points with
  the same imaginary value in $K$, since in
  Figure~\ref{fig:subdomain_in_H} they correspond to the boundary
  points that can be connected by a ray from the origin.
 \end{proof}

 \begin{lem}\label{lem:K}
  $\Stab(\bP^{1})/(\bZ)\bC$ is conformally equivalent to $\bC^{*}$.
 \end{lem}
 \begin{proof}
  This follows from Riemann mapping theorem and Reflection principle.
  We use notation from Figure \ref{fig:K}; the origin and the infinite
  point are called $A$ and $B$, the boundary lines $e^{\pm x}\cos y=1$
  are called $L_{\pm}$, and the upper imaginary axis by $L$.
  
  By $z\mapsto \frac{z-i}{z+i}$, $K$ is conformally equivalent to the
  subdomain in the unit disk ($|z|<1$) bounded by $L_{\pm}$.  Here, two
  boundary points of $K$ with the same imaginary value go to the
  boundary points with the same real value.  Let us denote the
  upper-half and lower-half of $K$ by $K_{u}$ and $K_{l}$.  See Figure
  \ref{fig:K_unit}.

  Next, by the Riemann mapping theorem there is a bijective conformal
  mapping from $K_{u}$ to the unit disk.  We can extend any isomorphism
  of bounded domains to a homeomorphism on their closures by
  \cite[Theorem 11-1]{COH}. By a linear fractional transformation, we
  can rearrange three points on the boundary in arbitrary way as long as
  we keep their order.  Hence $K_{u}$ is conformally equivalent to the
  unit disk where $A$ and $B$ correspond to $-1$, $1$, and the
  upper-half circle and the lower-half circle correspond $L_{-}$ and
  $L$.  See Figure \ref{fig:unit}.
  \setlength{\minitwocolumn}{0.45\textwidth}
  \addtolength{\minitwocolumn}{0.2\columnsep}
  \begin{figure}[H]
   \begin{tabular}{cc}
    \begin{minipage}{\minitwocolumn}
     \begin{center}
      \scalebox{.7}{
      \input{K_66.pstex_t}
      }
      \caption{$K$ in the unit disk}
      \label{fig:K_unit}
     \end{center}
    \end{minipage}
    &
    \begin{minipage}{\minitwocolumn}
     \begin{center}
    \scalebox{.7}{
      \input{K_8.pstex_t}
      }
      \caption{$K_{u}$ as the unit disk}
      \label{fig:unit}
     \end{center}
    \end{minipage}
   \end{tabular}
  \end{figure}
  
  Next, by the composition of $z\mapsto
  -2i\left(\frac{zi+1}{zi-1}\right)$, $z\mapsto
  \frac{z+\sqrt{z^{2}-4}}{z}$, and $z\mapsto \frac{1}{z}$, $K_{u}$ as
  the unit disk is conformally equivalent to the lower-half disk, where
  $B$ and $A$ correspond $-1$ and $1$, and the lower-half circle
  corresponds to $L_{-}$.  See Figure \ref{fig:lower}.
  \setlength{\minitwocolumn}{0.45\textwidth}
  \addtolength{\minitwocolumn}{0.3\columnsep}
  \begin{figure}[H]
   \begin{tabular}{c c}
    \begin{minipage}{\minitwocolumn}
     \begin{center}
      \scalebox{.7}{
    \input{K_11.pstex_t}
    }
    \caption{$K_{u}$ as the half disk}
    \label{fig:lower}
     \end{center}
    \end{minipage}
    &
    \begin{minipage}{\minitwocolumn}
     \begin{center}
    \scalebox{.7}{
    \input{K_12.pstex_t}
    }   
    \caption{$K$ as the unit disk}
     \end{center}
    \end{minipage}
   \end{tabular}
  \end{figure}

  Now, by the Reflection Principle we can extend the bijective
  conformal mapping from $K$ to the unit disk, where the two points on
  $L_{\pm}$ with the same imaginary part go to the points on the
  boundary on the unit disk with the same real part.

  Finally, $z\mapsto -i(\frac{z+1}{z-1})$ sends $K$ as the unit disk to
  the upper-half plane, where two points on the boundary with the same
  real value are mapped to two points on the boundary with the same
  absolute value. Then, $z \mapsto z^{2}$ sends $K$ as the upper-half
  plane to $\bC^{*}$ and identifies on the real axis.
 \end{proof}

  \def\proofname{Proof of Theorem 1.1}
  \begin{proof}
   By Lemma \ref{lem:K}, we have $\Stab(\D(\bP^{1}))/(\bZ)\bC\cong
   \bC^{*}$.  The action of $\bZ$ on
   $\mathfrak{X}=\Stab(\D(\bP^{1}))/\bC$ gives an exact sequence $0\to
   \pi_{1}(\mathfrak{X})\to \pi_{1}
   (\mathfrak{X}/\bZ)\stackrel{\alpha}{\to} \pi_{0}(\bZ) \to
   \pi_{0}(\mathfrak{X})$.  We will show that $\mathfrak{X}$ is
   connected.  Recall that $\Stab(\D(\bP^{1}))=(\bZ)\bC \cdot X$ and, by
   Lemma \ref{lem:1}, $X\cong \bH$ is connected, hence so is $\bC\cdot
   X$.  It remains to check that $(\bZ)$ fixes some connected component
   of $\Stab(\D(\bP^{1}))$, but it fixes $\{(Z, \cP)\in
   \Stab(\D(\bP^{1}))\mid \cP((0,1])=\Coh \bP^{1}\}$, which lives in
   $X$.  So $\mathfrak{X}$ is connected and the map $\alpha$ is a
   surjective map $\bZ\to \bZ$, therefore $\alpha$ is injective and
   $\pi_{1}(\mathfrak{X})=0$.  Hence $\mathfrak{X}$ is the universal
   covering of $\bC^{*}$; i.e., $\Stab(\D(\bP^{1}))/\bC\cong \bC$.
   Moreover, since $H^{1}(\bC, \cO)=0$, $\Stab(\D(\bP^{1}))\cong
   \bC^{2}$.
  \end{proof}
  \def\proofname{Proof}

  \section{Walls and hearts of $\Stab(\D(\bP^{1}))$}\label{sec:walls-hearts-stabdbp}
  We define the ``{\it cell}'' $\Stab_{\cC}(\cT)$ for a heart $\cC$ by
  $\Stab_{\cC}(\cT)=\{(Z,\cP)\in \Stab(\cT)\mid \cP((0, 1])=\cC\}$
  ($S_{\cC}(\cT)$ for short).  In this section, we describe all cells
  $S_{\bC}(\D(\bP^{1}))$ and how they fit together. We describe a
  fundamental domain of $\Stab(\D(\bP^{1}))/(\bZ)\bR$ as a real
  manifold, this shows that the quotient $\Stab(\D(\bP^{1}))/(\bZ)\bR$
  is an open torus.
  
  For the hearts $\cC_{j}$ and $\cC_{p, i, j}$ from Corollary
  \ref{cor:all_hearts}, let $S_{j}=S_{C_{j}}$ and
  $S_{p,i,j}=S_{\cC_{p,i,j}}$. Then they can be described as
  \begin{align*}
   S_{j}
   &=\{(Z, \cP)\in \Stab(\D(\bP^{1}))\mid \phi(\cO_{x}[j])=1, \   0<\phi(\cO[j])<1\}\cong \bR_{<0}\times \bH;\\
   S_{p, i, j}
   &=\{(Z, \cP)\in \Stab(\D(\bP^{1}))
   \mid 0< \phi(\cO(i-1)[p+j]), \phi(\cO(i)[j]) \leq 1\}
   \cong H^{2}.
  \end{align*}
  For $p\geq 1$ and $i,j\in \bZ$, we have the following 4-dimensional
  manifolds with boundaries:
  \begin{align*}
   S^{-}_{p, i, j}
   &=\{(Z, \cP)\in S_{p, i, j}\mid
   \phi(\cO(i-1)[p+j])\geq \phi(\cO(i)[j])\};\\
   S^{+}_{p, i, j}
   &=\{(Z, \cP)\in S_{p, i, j}\mid
    \phi(\cO(i-1)[p+j])\leq \phi(\cO(i)[j])\};\\
   P_{i, j}
   &=S^{+}_{1, i, j}
   \cup \left( \cup_{p>0} S_{p, i, j} \right).
  \end{align*}
  For each $S^{\pm}_{p,i,j}$, we will find which ones are {\it
  neighbors} in the sense that the intersection of their closures is a
  codimension-one submanifold.  See the following graphic descriptions
  of $S^{-}_{p, i, j}$ and $S^{+}_{p, i, j}$, and see Figure
  \ref{fig:std} for $S_{0}$.  \setlength{\minitwocolumn}{0.45\textwidth}
  \addtolength{\minitwocolumn}{0.3\columnsep}
  \begin{figure}[H]
   \begin{tabular}{c c}
    \begin{minipage}{\minitwocolumn}
     \begin{center}
      \scalebox{.7}{ 
      \input{K_13.pstex_t}
      }   
      \caption{$S^{-}_{p, i, j}$}
     \end{center}
    \end{minipage}
    &
    \begin{minipage}{\minitwocolumn}
     \begin{center}
     \scalebox{.7}{ 
      \input{K_14.pstex_t}
      }   
      \caption{$S^{+}_{p, i, j}$}
     \end{center}
    \end{minipage}
   \end{tabular}
  \end{figure}
  When $p=1$, let us simplify notations as follows:
  \begin{align*}
   S^{-}_{i, j}
   &=S^{-}_{1, i, j}; \  S^{+}_{i, j}=S^{+}_{1, i, j};  \ W_{i, j}=S^{-}_{i, j}\cap S^{+}_{i, j};\\
   l_{i, j}
   &=\{(Z, \cP)\in W_{i, j}\mid \phi(\cO(i-1)[j+1])=\phi(\cO(i)[j])=1\}.
  \end{align*}
  \begin{prop}\label{prop:walls}
   $W_{i, j}$ is a wall for any $i, j\in \bZ$,
   and there is no other walls.
  \end{prop}
  \begin{proof}
   In $S^{-}_{i, j}$, all line bundles and torsion sheaves are
   semistable, but in $S^{+}_{i, j}$, only two line bundles are
   semistable.  In the fundamental domain of $\Stab(\D(\bP^{1}))$ in
   Lemma \ref{lem:3}, we have only one wall $y=\pi$ that is the quotient
   of the set consisting of $W_{i, j}$ for all $i, j$.
  \end{proof}
  
  Before the neighbors lemma for $S^{-}_{i,j}$, we need the following
  technical corollaries from Proposition \ref{prop:exp}.

  \begin{cor}\label{cor:consecutive}
   Let $(Z, \cP)$ be a stability condition such that $\cO(i),
   \cO(i-1)\in \cP(\alpha)$ for some $i$ and $\alpha$, then all line
   bundles and torsion sheaves are semistable.
  \end{cor}
  \begin{proof}
   It follows from the case 3 in Proposition \ref{prop:exp}.
  \end{proof}

  \begin{cor}\label{cor:consecutive_2}
   There is no stability condition $(Z, \cP)$ such that
   $m(\cO(i)[j])=m(\cO(i-1)[j+1])(1+1/n)$ for some $n\in \bZ
   \cup \{\infty\}$ and $\phi(\cO(i)[j])=\phi(\cO(i-1)[j+1])-1$.
  \end{cor}
  \begin{proof}
   The assumption implies $Z(\cO_{x}[j])=(-1/n) Z(\cO(i-1)[j+1])$ since
   $Z(\cO(i)[j])+Z(\cO(i-1)[j+1])=Z(\cO_{x}[j])$.  Then
   $Z(\cO(i-1-n)[j+1])=0$ for $n\neq 0, \infty$, and $Z(\cO_{x})=0$ or
   $\infty$ for otherwise.  However, $\cO(i)[j+1], \cO(i-1)[j+1]\in
   \cP(\phi(\cO(i-1)[j+1]))$ implies all line bundles and torsion
   sheaves are semistable by Corollary \ref{cor:consecutive}.
  \end{proof}

   \begin{lem}\label{lem:neighbor_1}
    $S^{-}_{i, j}$ has neighbors: (1) $S^{-}_{i+1,j}$; (2)
    $S^{-}_{i-1,j}$; (3) $S^{+}_{i,j}$.  Moreover, in each case we can
    describe the intersection of their closures as a submanifold with
    boundaries: (1) $\{(Z, \cP)\in S^{-}_{i+1,j}\mid
    \phi(\cO(i)[j])=0\}$ with boundaries $l_{i+n, j}$ for $n\geq 0$ and
    $l_{i-n, j+1}$ for $n>0$; (2) $\{(Z, \cP)\in S^{-}_{i,j}\mid
    \phi(\cO(i-1)[j+1])=1\}$ with boundaries $l_{i+n, j}$ for $n>0$ and
    $l_{i-n, j+1}$ for $n\geq 0$; (3) $W_{i,j}$ with boundaries $l_{i,
    j}$ and $l_{i, j+1}$.  See the following figure of $S^{-}_{i, j}$
    for a fixed $m(\cO(i-1)[j])$.
  \begin{figure}[H]
   \begin{center}
    \scalebox{.7}{
    \input{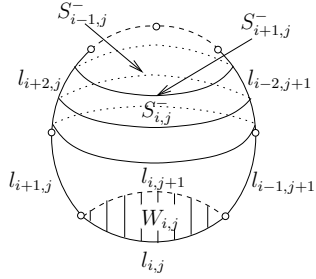}
    }   
    \caption{$S^{-}_{i, j}$ is the
    inside of the figure above.}
   \end{center}
  \end{figure} 
   \end{lem}

  \begin{proof}
   Let $\sigma=(Z, \cP)\in S^{-}_{i, j}$.  If we rotate $Z(\cO(i)[j])$
   counterclockwise, $\sigma$ will be in $W_{i, j}$ when
   $\phi(\cO(i)[j])=\phi(\cO(i-1)[j+1])$.
   
   If $\phi(\cO(i-1)[j+1])\neq 1$, we see that $\sigma$ will be in
   $S^{-}_{i+1, j}$ when we rotate $Z(\cO(i)[j])$ clockwise, and $\sigma$
   will be in $S^{-}_{i-1, j}$ when we rotate $Z(\cO(i-1)[j+1])$
   counterclockwise.
  \setlength{\minitwocolumn}{0.45\textwidth}
  \addtolength{\minitwocolumn}{0.3\columnsep}
  \begin{figure}[H]
   \begin{tabular}{c c}
    \begin{minipage}{\minitwocolumn}
     \begin{center}
      \scalebox{.7}{ 
      \input{K_23.pstex_t}
      }   
      \caption{$S^{-}_{i+1, j}$}
     \end{center}
    \end{minipage}
    &
    \begin{minipage}{\minitwocolumn}
     \begin{center}
     \scalebox{.7}{ 
      \input{K_24.pstex_t}
      }   
      \caption{$S^{-}_{i-1, j}$}
     \end{center}
    \end{minipage}
   \end{tabular}
  \end{figure}
   Consider the case $\phi(\cO(i-1)[j+1])=1$.  Let us see what happens
   when we have $\phi(\cO(j)[j])=0$ by rotating $Z(\cO(i)[j])$
   clockwise. We can assume $m(\cO(i)[j])\neq m(\cO(i-1)[j+1])(1+ 1/n)$
   for any $n\in \bZ$ by Corollary \ref{cor:consecutive_2}.  Let us fix
   $m(\cO(i-1)[j+1])=1$, $i=j=0$.  If $1-1/(n-1)<m(\cO)< 1-1/n$ for some
   $n>1$, $m(\cO_{x})$ will be $1-m(\cO)$ and $\phi(\cO_{x})$ will be
   one.  Moreover, $\phi(\cO(-1+(n-1))[1])$ and $\phi(\cO(n-1))$ will be
   one, hence $\sigma$ will be in $l_{n-1, 0}$.  For general
   $m(\cO(i-1)[j+1])$, $i$ and $j$, $\sigma$ will be in $l_{n-1+i, j}$.
   On the other hand, if $1+1/n< m(\cO)< 1+1/(n-1)$ for some $n>0$, then
   by the same manner, $\sigma$ will be in $l_{-n, 1}$.  For general
   $m(\cO(i-1)[j+1])$, $i$ and $j$, $\sigma$ will be in $l_{-n+i, j+1}$.
   Some examples when $\phi(\cO)$ becomes zero can be seen in the
   following figures.  
   \setlength{\minitwocolumn}{0.45\textwidth}
   \addtolength{\minitwocolumn}{0.3\columnsep}
   \begin{figure}[H]
    \begin{tabular}{c c}
     \begin{minipage}{\minitwocolumn}
      \begin{center}
       \scalebox{.7}{ 
       \input{K_25.pstex_t}
       }   
       \caption{$l_{1, 0}$}
      \end{center}
     \end{minipage}
     &
     \begin{minipage}{\minitwocolumn}
      \begin{center}
       \scalebox{.7}{ 
       \input{K_26.pstex_t}
       }   
       \caption{$l_{-1, 1}$}
      \end{center}
     \end{minipage}
    \end{tabular}
   \end{figure}
  \end{proof}
  \begin{lem}\label{lem:neighbor_2}
   For $p>1$, $S^{-}_{p, i, j}$ has neighbors $S^{+}_{p,i,j}$,
   $S^{+}_{p-1,i,j+1}$, and $S^{+}_{p-1, i, j}$.  For $p>0$, $S^{+}_{p, i,
   j}$ has neighbors $S^{-}_{p, i, j}$, $S^{-}_{p+1, i, j-1}$,
   and  $S^{-}_{p+1, i, j}$.  See the figure below for some fixed
   $m(\cO(i-1)[j])$ and $m(\cO(i)[j+1])$.
  \begin{figure}[H]
   \begin{center}
    \scalebox{.7}{
    \input{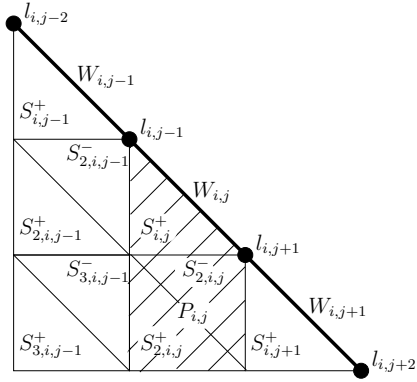}
    }   
    \caption{Some components connected to $S^{+}_{i, j}$}
   \end{center}
  \end{figure} 
  \end{lem}
  \begin{proof}
   Let $\sigma=(Z, \cP)\in S^{-}_{p, i, j}$. If we rotate $Z(\cO(i)[j])$
   counterclockwise, then $\sigma$ will be in $S^{+}_{p, i, j}$ when
   $\phi(\cO(i)[j])=\phi(\cO(i-1)[p+j])$.  If we rotate $Z(\cO(i)[j])$
   clockwise, $\sigma$ will be in $S^{+}_{p-1, i, j+1}$ when
   $\phi(\cO(i)[j])=0$; i.e., $\phi(\cO(i)[j+1])=1$.  If we rotate
   $Z(\cO(i-1)[p+j])$ counterclockwise, $\sigma$ will be in $S^{+}_{p-1,
   i, j}$ when $\phi(\cO(i-1)[p+j])>1$; i.e.,  $1\geq
   \phi(\cO(i-1)[p-1+j])>0$.

   Let $\sigma=(Z, \cP)\in S^{+}_{p, i, j}$ for $p>0$. By the same manner
   above, $\sigma$ will be in $S^{-}_{p+1, i, j-1}$ by rotating
   $Z(\cO(i)[j])$ counterclockwise, and $\sigma$ will be in
   $S^{-}_{p+1, i, j}$ by rotating $Z(\cO(i-1)[p+j])$ clockwise.
  \end{proof}
  \begin{prop}
   A fundamental domain of $\Stab(\D(\bP^{1}))/(\bZ) \bR$ is described
   by Figure \ref{fig:ball}.  In particular
   $\Stab(\D(\bP^{1}))/(\bZ)\bR$ is isomorphic to an open torus.
  \end{prop}
  \begin{proof}
   Let us see how $S^{-}_{n, j}$ and $P(n, j)$ for all $n\in \bZ$, and
   $S_{j}$ and $S_{j+1}$ fit together.  Let $B_{i, j}$ be the union of
   them.  We attach $S^{-}_{i+1, j}$ to the front-side of $S^{-}_{i,
   j}$, $S^{-}_{i-1, j}$ to the backside of $S^{-}_{i, j}$, and $P(i,
   j)$ to $W_{i, j}$.  Then we get Figure \ref{fig:B}, where we attach
   $S_{j}$ to the front-side and $S_{j+1}$ to the backside.  Moreover,
   it is isomorphic to an open ball in Figure \ref{fig:ball}.
   Therefore, once we identify $S_{j}$ and $S_{j+1}$, we see that
   $\Stab(\D(\bP^{1}))/(\bZ)\bR$ is an open torus.
   \setlength{\minitwocolumn}{0.45\textwidth}
   \addtolength{\minitwocolumn}{0.3\columnsep}
   \begin{figure}[H]
    \begin{tabular}{c l}
     \begin{minipage}{\minitwocolumn}
      \begin{center}
       \scalebox{.7}{ 
       \input{K_188.pstex_t}
       }   
       \caption{$B_{i, j}$}
       \label{fig:B}
      \end{center}
     \end{minipage}
     &
     \begin{minipage}{\minitwocolumn}
      \begin{center}
       \scalebox{.7}{ 
       \input{K_200.pstex_t}
       }   
       \caption{$B_{i, j}$ as a ball}
       \label{fig:ball}
       \end{center}
     \end{minipage}
    \end{tabular}
   \end{figure}
  \end{proof}    
  \begin{rmk}
   Heart $\Coh \bP^{1}=\cC_{0}$ and its shifts $\cC_{j}$ are the most
   special hearts that appear in $\Stab(\D(\bP^{1}))$, in the sense that
   the dimension of $S_{j}$ is the smallest. They are also the most
   symmetric ones; i.e., $\cC_{j}$ is fixed by $\mathrm{Pic}(\bP^{1})$
   and $\Aut(\bP^{1})$, while any other heart is fixed only by
   $\Aut(\bP^{1})$.  The way all $S_{\cC}$ fit into $\Stab(\D(\bP^{1}))$
   gives a picture of degenerations of the hearts $\cC_{j}$ and
   $\cC_{i,j}$, that are of the homological dimension one, to the hearts
   $\cC_{p,i,j}$ for $p>1$, that are of the homological dimension zero.
  \end{rmk}
 
\end{document}